\documentclass[11pt]{article}
\usepackage{natbib}
\usepackage{epsfig}       
\begin{document}
\begin{center}
{\LARGE {\bf Towards Reconciliation between Bayesian  \\

\vspace{0.3cm}

and Frequentist Reasoning}}\\
\end{center}
\vspace{1cm}

\begin{center}
Toma\v z Podobnik$^{1,\;2} $ {\small and} 
 Tomi \v Zivko$^2$ {\small (presenter)}\\
E-mail: tomaz.podobnik@ijs.si \, \, tomi.zivko@ijs.si \\
$^1$ {\it Physics Department, University of Ljubljana, Slovenia}\\
$^2$ {\it ''Jo\v zef Stefan'' Institute, Ljubljana, Slovenia}
\end{center}
\begin{abstract}

A theory of quantitative inference about the parameters of sampling
distributions is constructed deductively by following
very general rules, referred to as the Cox-P\'{o}lya-Jaynes Desiderata.
The inferences are made in terms of probability distributions that 
are assigned to the parameters. The Desiderata, focusing primarily on
consistency of the plausible reasoning, lead to unique assignments
of these probabilities in the case of sampling distributions that
are invariant under Lie groups. In the scalar cases, e.g. in the
case of inferring a single location or scale parameter, the
requirement for logical consistency is equivalent to the requirement
for calibration: the consistent probability distributions are
automatically also the ones with the exact calibration and
\textit{vice versa}. This equivalence speaks in favour of 
reconciliation between the Bayesian and the frequentist schools 
of reasoning.

\end{abstract}
\section{Introduction}

A theory of quantitative inference about the parameters of sampling
distributions is formulated with special attention being paid
to the consistency of the theory and to its ability to make verifiable
predictions. In the present article only basic concepts of the 
theory and their most important applications are presented while
details can be found elsewhere~[1].

Let $p(x_1|\theta I)$ be the probability for a random variate $x$ to 
take the value $x_1$ (to take a value in an interval $(x_1,x_1+dx)$ in
the case of a continuous variate), given the family $I$ of sampling 
distributions, and the value $\theta$ of the parameter that specifies 
a unique distribution within the family (for example, a sampling 
distribution from the exponential family $I$, 
$\tau^{-1}\exp{\{-x/\tau \}}$, is uniquely determined by the value 
of the parameter $\tau$).
An inference about the parameter is made by specifying a real number,
called (degree of) plausibility, $(\theta|x_1x_2\ldots I)$, to 
represent our degree of belief in the value of the (continuous) parameter 
to be within an interval 
$(\theta,\theta+d\theta)$. Every such plausibility is conditioned upon 
the information that consists of measured value(s) 
$x_1$, $x_2$, $\ldots$ of the sampling variate and of the 
specified family $I$ of sampling distributions.

We assume all considered plausibilities to be 
subjects to very general requirements, referred to as the 
Cox-P\'{o}lya-Jaynes (CPJ) Desiderata~[2,1], focusing
primarily on consistency of the plausible reasoning.
The requirement of consistency can be regarded as the first of the
requirements to be satisfied by every theoretical system,
be it empirical or non-empirical. As for an empirical
system, however, besides being consistent, it must also be  
falsifiable~[3].
We therefore added a Desideratum to CPJ Desiderata, requiring that
the predictions of the theory must be verifiable so that, in
principle, they may be refuted.

It should be stressed that in this way the list of basic rules is
completed. That is, the entire theory of inference about the
parameters is built deductively from the aforementioned Desiderata:
in order not to jeopardize the consistency of the theory no additional
\textit{ad hoc} principles are invoked. 

\section{Cox's and Bayes' Theorems}

Richard Cox showed~[4] that a system for manipulating
plausibilities is either isomorphic to the probability system or
inconsistent (i.e. in contradiction with CPJ Desiderata). Without
any loss of generality we therefore once and for all choose 
probabilities $p(\theta|x_1I)$ among all possible plausibility
functions $(\theta|x_1I)$ to represent our degree of belief in 
particular values of inferred parameters. In this way the so-called inverse
probabilities, $p(\theta|x_1I)$, and the so-called direct
(or sampling) probabilities $p(x_1|\theta I)$, become 
subjects to identical rules.

Transformations of probability distributions that are induced
by variate transformations are also uniquely determined
by the Desiderata. Let $f(x|\theta I)$ be the probability density 
function (pdf) for a continuous random variate $x$ so that its
probability distribution is expressible as
\begin{equation}
 p(x|\theta I)=f(x|\theta I)\,dx \ .
\end{equation}
Then, if the variate $x$ is subject to a one-to-one transformation
$x\longrightarrow y=g(x)$, the pdf for $y$ reads:
\begin{equation}
 f(y|\theta I')=f(x|\theta I)\,\Bigl|\frac{dy}{dx}\Bigr|^{-1} 
\label{eq:vartr1}
\end{equation}
(by using the symbol $I'$ instead of $I$ on the left-hand side of
(\ref{eq:vartr1}) it is stressed that the above transformations 
may in general alter the form of the sampling distribution).   
Since the direct and the inverse probabilities are 
subjects to the same rules, the transformation of the pdf
for the inferred parameter, 
$f(\theta|x I)$, under a one-to-one transformation 
$\theta\longrightarrow \nu = \bar{g}(\theta)$ is analogous to the
transformation of the sampling pdf:
\begin{equation}
f(\nu | x I) = f(\theta | x I)\,
\Bigl|\frac{d\nu}{d\theta}\Bigr|^{-1} \ .
\label{eq:vartr2}
\end{equation}

Once the probabilities are chosen,
the usual product and sum rules~[2] become the fundamental
equations for manipulating the probabilities, while many other equations 
follow from the repeated applications of the two. 
In this way, for example, Bayes' Theorem for updating the
probabilities can be obtained: 
\begin{equation}
f(\theta | x_1 x_2 I ) = 
\frac{f(\theta |x_1 I)\,p(x_2 | \theta x_1 I)}
     {\int f(\theta' |x_1 I)\,p(x_2 | \theta' x_1 I)\,d\theta'} \ .
\label{eq:bayesth}
\end{equation}
Here $f(\theta |x_1 I)$ denotes the pdf for $\theta$ based on $x_1$ 
and $I$ only (i.e. prior to taking datum $x_2$ into account), 
$p(x_2 | \theta x_1 I)$ is the 
probability for $x_2$ (the so-called likelihood) given values
$\theta$ and $x_1$, while the integral in the denominator on 
the right-hand side ensures appropriate normalization of the updated 
pdf $f(\theta | x_1 x_2 I )$ for $\theta$ (i.e. the 
pdf for $\theta$ posterior to taking $x_2$ into account).

Bayes' Theorem (\ref{eq:bayesth}) allows only for updating 
pdf's $f(\theta|x_1I)$ that were already assigned prior to their updating. 
Consequently, the existing applications of our basic rules must
be extended in order to allow for assignment of probability
distributions to the parameters, with such assignments representing 
natural and indispensable starting points in every sequential
updating of probability distributions.

\section{Consistency Theorem}

According to the CPJ Desiderata, the pdf for $\theta$
should be invariant under reversing the order of taking into account 
two independent measurements of the sampling variate $x$. This is 
true if and only if the pdf that is assigned to $\theta$ on the
basis of a single measurement of $x$, is directly proportional
to the likelihood for that measurement,
\begin{equation}
 f(\theta |x I) \; =
 \frac {\pi(\theta)\,p(x | \theta I)}
       {\int\pi(\theta')\,p(x | \theta' I)\,d\theta'} \ ,
 \label{eq:consth}
\end{equation}
where $\pi(\theta)$ is the consistency factor while the integral in
the denominator on the right-hand side of (\ref{eq:consth}) again
ensures correct normalization of $f(\theta |x I)$.

There is a remarkable similarity between the Bayes' Theorem 
(\ref{eq:bayesth}), applicable for \textit{updating} the probabilities, 
and the Consistency Theorem (\ref{eq:consth}), applicable
for \textit{assigning} the probability distributions to the 
values of the inferred parameters, but there is also a fundamental
and very important difference between the two. While $f(\theta|x_1I)$
in the former represents the pdf for $\theta$ prior to taking
datum $x_2$ into account, $\pi(\theta)$ in the latter is (by 
construction of the Consistency Theorem~[1]) 
just a proportionality coefficient 
between the pdf for $\theta$ and the appropriate likelihood
$p(x|\theta I)$, so that no probabilistic inference is 
ever to be made on the consistency factor alone, nor can $\pi(\theta)$
be subject to the normalization requirement that is otherwise
perfectly legitimate in the case of prior pdf's.

The form of the consistency factor depends on the only relevant
information that we posses before the first datum is collected, i.e.
it depends on the specified sampling model. Consequently, when
assigning probability distributions to the parameters of 
the sampling distributions from the same family $I$, this must be
made according to the Consistency Theorem by using
the consistency factors of the forms that are identical
up to (irrelevant) multiplication constants.

\section{Consistency Factor}

According to (\ref{eq:vartr2}) and (\ref{eq:consth}) combined, the
consistency factors $\pi(\theta)$ for $\theta$ and 
$\widetilde{\pi}\bigl(\bar{g}(\theta)\bigr)$ for the
transformed parameter $\bar{g}(\theta)$ are related as
\begin{equation}
 \widetilde{\pi}\bigl(\bar{g}(\theta)\bigr)
 = k\,\pi(\theta)\,\bigl|\bar{g}'(\theta)\bigr|^{-1} \ , 
\label{eq:feq1}
\end{equation}
where $k$ is an arbitrary constant (i.e. its value is independent
of either $x$ or $\theta$), while $\bar{g}'(\theta)$
denotes the derivative of $\bar{g}(\theta)$ with respect to
$\theta$. However, for the parameters of sampling distributions 
with the form $I$ that is invariant under simultaneous transformations 
$g_a(x)$ and $\bar{g}_a(\theta)$ of the sample and the parameter
space, 
\begin{equation}
 f\bigl(g_a(x)|\bar{g}_a(\theta) I'\bigr)
 \!=\!\!
 f(x|\theta I)\,\bigl|g_a\!'(x)\bigr|^{-1} 
 \!\!\!=\!  
 f\bigl(g_a(x)|\bar{g}_a(\theta) I\bigr)
\end{equation}
(i.e. when $I'=I$), $\widetilde{\pi}$ and $\pi$ must be identical
functions up to a multiplication constant, so that (\ref{eq:feq1})
reads:
\begin{equation}
 \pi\bigl(\bar{g}_a(\theta)\bigr)
 = k(a)\,\pi(\theta)\,\bigl|\bar{g}_a'(\theta)\bigr|^{-1} \ . 
\label{eq:feq2}
\end{equation}
Index $a$ in the above expressions indicates parameters of the
transformations and $k$, in general, can be a function of $a$.  
In the case of multi-parametric transformation groups the derivative
$\bar{g}_a'(\theta)$ is to be substituted by the appropriate Jacobian.

The above functional equation has a unique solution for the
transformations $\bar{g}_a(\theta)$ with the continuous 
range of admissible values $a$, i.e. if the set of
admissible transformations $\bar{g}_a(\theta)$ forms
a Lie group. If a sampling distribution for $x$ is invariant
under a Lie group, then it is necessarily reducible (by
separate one-to-one transformations of the sampling variate $x \to y$
and of the parameter $\theta \to \mu$) to a sampling
distribution that can be expressed as a function of 
a single variable $y-\mu$, $f(y|\mu I) = \phi(y-\mu)$.
Sampling distributions of the form $\sigma^{-1}\psi(x/\sigma)$
are examples of such distributions: by substitutions
$y=\ln{x}$ and $\mu=\ln{\sigma}$ they transform into
$\phi(y-\mu)=\exp{\{y-\mu\}}\,\psi\bigl(\exp{\{y-\mu\}}\bigr)$
(the scale parameters $\sigma$ are reduced to location parameters
$\mu$). 

It is therefore sufficient to determine
the form of consistency factors for the location parameter $\mu$
since we can always make use of (\ref{eq:feq1}) to 
transform $\widetilde{\pi}\bigl(\mu=\bar{g}(\theta)\bigr)$ 
into the appropriate consistency factor $\pi(\theta)$
for the original parameter $\theta$.
Sampling distributions of the form $\phi(x-\mu)$ are
invariant under simultaneous translations 
$x\to x+a$ and $\mu\to\mu+a$; $\forall a\in(-\infty,\infty)$,
and the functional equation (\ref{eq:feq2}) in the case of the
translation group reads
\begin{equation}
 \pi(\mu+a) = k(a)\,\pi(\mu) \ ,
\end{equation}
implying the consistency factor for the location parameters to be 
$\pi(\mu)\propto\exp{\{-q\mu\}}$, with $q$ being an arbitrary 
constant. Accordingly, $\pi(\sigma)\propto\sigma^{-(q+1)}$
is the appropriate form of the consistency factor for the scale
parameters.

The value of $q$ is then uniquely determined by recognizing the
fact that sampling distributions of the forms $\phi(x-\mu)$ 
and $\sigma^{-1}\psi(x/\sigma)$ are just special cases of
two-parametric sampling distributions
\begin{equation}
 f(x|\mu\sigma
 I)=\frac{1}{\sigma}\,\phi\Bigl(\frac{x-\mu}{\sigma}\Bigr) \ ,
\label{eq:locsc}
\end{equation}
with $\sigma$ being fixed to unity and with $\mu$ being fixed to
zero, respectively. The consistency factor $\pi(\mu)$ therefore
corresponds to assigning pdf's $f(\mu|\sigma xI)$ while $\pi(\sigma)$
is to be used when assigning $f(\sigma|\mu xI)$. When neither $\sigma$ nor
$\mu$ is fixed, however, the pdf
(\ref{eq:locsc}) is invariant under a two-parametric group of 
transformations, $x\to ax+b$, $\mu\to a\mu+b$ and $\sigma\to a\sigma$;
$\forall a\in(0,\infty)$ and $\forall b\in(-\infty,\infty)$, and the 
functional equation (\ref{eq:feq2}) for the consistency factor 
$\pi(\mu,\sigma)$ for assigning $f(\mu\sigma|xI)$ reads
\begin{equation}
 \pi(a\mu+b,a\sigma)=\frac{k(a,b)}{a^2}\,\pi(\mu,\sigma) \ ,
\end{equation}
so that $\pi(\mu,\sigma)$ is to be proportional to $\sigma^{-r}$, $r$
being an arbitrary constant. According to the product rule,
$f(\mu\sigma|xI)$ can be factorized as
\begin{equation}
 f(\mu\sigma|xI)=f(\mu|\sigma xI)\,f(\sigma|xI) 
                =f(\sigma|\mu xI)\,f(\mu|xI) \ ,
\label{eq:product}
\end{equation}
where $f(\sigma|xI)$ and $f(\mu|xI)$ are the marginal pdf's, e.g.
\begin{equation}
 f(\sigma|xI)=\int f(\mu'\sigma|xI)\,d\mu' \ .
\label{eq:marginal} 
\end{equation}
The equalities (\ref{eq:product}) are achieved if and only if 
$q=0$ and $r=1$, i.e. if the three consistency factors, 
determined uniquely up to arbitrary multiplication constants, read:
\begin{equation}
 \pi(\mu) = 1 \ \ \hbox{and} \ \ 
 \pi(\sigma) = \pi(\mu,\sigma) = \sigma^{-1} \ .
\label{eq:cf1}
\end{equation} 

\section{Calibration}

In order to exceed the level of a mere speculation, the theory of
probabilistic inference about the parameters must be able to make 
predictions that can be verified (or falsified) by experiments.
Let therefore a random variate $x$ be subject to a family of 
sampling distributions $I$ and let several independent values $x_i$
of the variate be recorded. The predictions of the theory are made
in probabilities
\begin{equation}
 P\bigl(\theta\in(\theta_{i,1},\theta_{i,2})|x_i I\bigr) =
 \!\int_{\theta_{i,1}}^{\theta_{i,2}}\!\! f(\theta'|x_i I)\,d\theta' = 
 \delta 
\label{eq:calib1}
\end{equation}
that given measured value $x_i$ of the sampling variate,
an interval $(\theta_{i,1},\theta_{i,2})$ contains the actual
value of the parameter $\theta$ of the sampling distribution. For the
sake of simplicity, the intervals are chosen in such a way that
the probabilities $\delta$ are equal in each of the assignments.
The predictions are then verifiable at long term relative
frequencies: our probability judgments (\ref{eq:calib1}) are said
to be calibrated if the fraction of inferences with the 
specified intervals containing the actual value of the parameter,
coincides with $\delta$.

An exact calibration of an inference about a parameter $\theta$ is
ensured if the assigned pdf $f(\theta|x I)$ is related
to the (cumulative) distribution function $F(x,\theta)$
of the sampling variate as
\begin{equation} 
 f(\theta|x I) = \biggl| \frac{\partial}{\partial \theta} F(x,\theta)
                 \biggr| \ ,
\label{eq:calib2}
\end{equation}
and the consistency factors $\pi(\mu)$ and $\pi(\sigma)$
(\ref{eq:cf1}) do meet the above requirement. Furthermore, 
if besides of being calibrated (\ref{eq:calib2}), the pdf for 
$\theta$ is to be assigned
according to the Consistency Theorem (\ref{eq:consth}), the distribution
of the sampling variate $x$ is necessarily reducible 
to a distribution of the form $\phi(y-\mu)$ ~[5]. But exactly the same
necessary condition was obtained by requiring invariance of the 
sampling distribution under a Lie group, with such an invariance
being indispensable for determination of consistency factors
solely by imposing consistency to the assignment of pdf's. Imposing logical
consistency to the theory is thus equivalent to imposing
calibration to its predictions: every probabilistic 
inference about a parameter of a sampling distribution
that we are sure is consistent
will thus at the same time also be calibrated and, \textit{vice versa},
every calibrated inference, based on a posterior pdf
that is factorized according to (\ref{eq:consth}), will 
simultaneously be logically consistent, too. The equivalence
of the two requirements speaks in favour of reconciliation
between the (objective) Bayesian and the frequentist schools
of reasoning, the former paying attention primarily to logical 
consistency and the latter stressing the importance of verifiable 
predictions.

\section{Consistency Lost and Regained}

Numerous examples can be found with the sampling distributions lacking
invariance under Lie groups: there are sampling distributions 
for continuous random variates (e.g. the Weibull distribution) that are 
not invariant under continuous
groups of transformations, the symmetry can be broken by imposing
constraints to parameter spaces of otherwise invariant sampling 
distributions, or the sampling space may be discrete (e.g. in counting 
experiments), just to name three of the most common ones. No
consistent qualitative parameter inference is possible in such cases,
but under very general conditions the remedy is just to collect more
data relevant to the estimated parameters. Then, according to
the Central Limit Theorem, the discrete sampling distributions 
approach their dense (Gaussian) limits, the constraints of the
parameter spaces become more and more irrelevant, and the
sampling distributions of the maximum likelihood estimates of the
inferred parameter $\theta$ gain Gaussian shapes with $\theta$
being the location parameters of the latter, so the ability of making 
consistent inferences is regained.

\section{Consistency Preserved}

Consistency factors are determined exclusively by utilizing the tools
such as the product rule (\ref{eq:product}) and marginalization 
(\ref{eq:marginal}), that are deducible directly from the basic
Desiderata: in order to preserve consistency of inference it is
crucial to refrain from using \textit{ad hoc} shortcuts on the 
course of inference. For regardless
how close to our intuitive reasoning these \textit{ad hoc}
procedures may be, how well they may have performed in some other previous
inferences, and how respectable their names may sound (e.g. the principle
of insufficient reason or its sophisticated version - the principle of maximum
entropy, the principle of group invariance, the principle of maximum
likelihood, and the principle of reduction),
they are all found in general to lead to inferences that are neither
consistent nor calibrated.

\vspace{-0.15mm}

\section*{Acknowledgments}
We thank prof. Louis Lyons for making it possible to present
our research at PHYSTAT05 and to prof.~Ale\v s Stanovnik for 
proofreading a preliminary draft of this report and correcting
several inaccuracies.

\end{document}